\documentclass{article}
\usepackage{amsmath}
\usepackage{amssymb}
\usepackage{amscd}

\def\be{\begin{equation}}       \def\ee{\end{equation}}
\def\bd{\begin{displaymath}}    \def\ed{\end{displaymath}}
\def\beq{\begin{eqnarray}}      \def\eeq{\end{eqnarray}}
\def\bseq{\begin{eqnarray*}}    \def\eseq{\end{eqnarray*}}
\def\ba{\begin{array}}          \def\ea{\end{array}}
\def\ben{\begin{enumerate}}     \def\een{\end{enumerate}}

\def\lmt{\longmapsto}

\def\cop{\Delta}

\def\cnt{\varepsilon}
\def\ot{\otimes}

\def\GLqtwo{{GL}_{q}(2)}
\def\GLhtwo{{GL}_{h}(2)}

\def\GLpqtwo{{GL}_{p,q}(2)}
\def\GLhh'two{{GL}_{h,h'}(2)}

\def\Grs{{G}_{r,s}}
\def\Gmk{{G}_{m,k}}
\def\Grss'{{G}_{r}^{s,s'}}
\def\Gmkk'{{G}_{m}^{k,k'}}
\def\GLone{{GL}({1})}
\def\GLtwo{{GL}({2})}

\def\rinv{{r}^{-1}}

\def\ident{{\bf 1}}

\begin{document}
\begin{center}
{\bf
ON SOME TWO PARAMETER QUANTUM AND JORDANIAN DEFORMATIONS, AND THEIR
COLOURED EXTENSIONS\footnote{Presented by D.P. at the LMS - Durham
Symposium on Quantum Groups, Durham, 19-29 July 1999, to appear in the
Proceedings.}
}

\bigskip
\bigskip
\bigskip

{\large\bf Deepak Parashar\footnote{email: deeps@scms.rgu.ac.uk} 
and Roger J. McDermott\footnote{email: rm@scms.rgu.ac.uk}
}

\bigskip
{\sl School of Computer and Mathematical Sciences\\
The Robert Gordon University, St. Andrew Street\\
Aberdeen AB25 1HG, United Kingdom\\
}
\bigskip
\bigskip
{\bf Abstract}
\end{center}
\medskip

This paper suveys some recent algebraic developments in two parameter
Quantum deformations and their Nonstandard (or Jordanian) counterparts. In
particular, we discuss the contraction procedure and the quantum group
homomorphisms associated to these deformations. The scheme is then set in
the wider context of the coloured extensions of these deformations,
namely, the so-called Coloured Quantum Groups.
\par

\section{Introduction}

Recent years have witnessed considerable development in the study of
multiparameter quantum deformations from both, the algebraic as well the
differential geometric point of view. These have also found profound
applications in many diverse areas of Mathematical Physics. Despite of the
intensive and successful development of the mathematical theory of
multiparameter quantum deformations or quantum groups, various important
aspects still need thorough investigation. Besides, all quantum groups
seem to have a natural coloured extension thereby defining corresponding
coloured quantum groups. It is the aim of this paper to address some of
the key issues involved.\\
Two parameter deformations provide an obvious step in constructing
generalisations of single parameter deformations. Besides being
mathematically interesting in their own right, two parameter quantum
groups serve as very good examples in generalising physical theories based
on the quantum group symmetry. $\GLpqtwo$ and $\GLhh'two$ are well known
examples of two parameter Quantum and Jordanian deformations of the space
of $2 \times 2$ matrices. Just as both these quantum groups are of great
significance in building up various mathematical and physical theories, it
is worthwhile to look for other possible examples, including the
`coloured' ones, which might play a fundamental role in future researches.
We wish to focus our attention on a new two parameter quantum group [1],
$\Grs$ which sheds light on some of the above mentioned issues. $\Grs$ is
a quasitriangular Hopf algebra generated by five elements, four of which
form a Hopf subalgebra ismomorphic to $\GLqtwo$, while the fifth generator
relates $\Grs$ to $\GLpqtwo$. \\
The $\Grs$ quantum group, which is the basis of our investigation, is
defined in Section II. In Section III, we give a new Jordanian analogue of
$\Grs$, denoted $\Gmk$ and establish a homomorphism with $\GLhh'two$. Both
$\Grs$ and $\Gmk$ admit a natural coloured extension and this is given in
Section IV. Section V generalises the contraction procedure to the case of
coloured quantum groups and discusses various homomorphisms. In section
VI, we make concluding remarks and give possible physical significance of
our results. Throughout this paper, we shall endeavour to refrain from too
much of technical details, which can be found in the appropriate
references.

\section{Two parameter $q$- deformations}

The quantum group $\Grs$ was defined in [1] as a quasitriangular Hopf
algebra with two deformation parameters $r$ and $s$, and generated by
five elements $a$,$b$,$c$,$d$ and $f$. The generators $a$,$b$,$c$,$d$ of
this Hopf algebra form a subalgebra, infact a Hopf subalgebra, which
coincides exactly with the single parameter dependent $\GLqtwo$ quantum
group when $q=r^{-1}$. Moreover, the two parameter dependent $\GLpqtwo$
can also be realised through the generators of this $\Grs$ Hopf algebra,
provided the sets of deformation parameters $(p,q)$ and $(r,s)$ are
related to each other in a particular fashion. This new algebra can,
therefore, be used to realise both $\GLqtwo$ and $\GLpqtwo$ quantum
groups. Alternatively, this $\Grs$ structure can be considered as a two
parameter quantisation of the classical $\GLtwo\ot \GLone$ group.  The
first four generators of $\Grs$, i.e. $a$, $b$, $c$, $d$ correspond to
$\GLtwo$ group at the classical level and the remaining generator $f$ is
related to $\GLone$. In fact, $\Grs$ can also be interpreted as a quotient
of multiparameter $q$- deformation of $GL(3)$.
\par
The elements of $\Grs$ can be conveniently arranged in the matrix 
$T=\left( \begin{smallmatrix}a&b&0\\c&d&0\\0&0&f\end{smallmatrix}
\right)$
and the coalgebra and counit are $\cop (T)=T\dot{\ot} T$, 
$\cnt (T)=\ident$. It should be mentioned that the quantum determinant 
$\delta = D f$ (where $D = ad-\rinv bc$) is group-like but not
central. The above block diagonal form of the $T$- matrix is particularly
convenient to understand the related schematics. The $\Grs$ $R$-matrix is
given in [1,2], and the most general Hopf algebra generated by this $R$-
matrix is multiparameter $GL(3)$ with the $T$-matrix of the form
$\left(
\begin{smallmatrix}a&b&x_{1}\\c&d&x_{2}\\v_{1}&v_{2}&f\end{smallmatrix}
\right)$.
It can be shown that the two-sided Hopf ideal generated by $x_{1},x_{2}$
when factored out yields the Inhomogenous multiparameter
$IGL(2)$. Furthermore, if one factors out yet another two-sided Hopf ideal
generated by elements $v_{1},v_{2}$, what one obtains is precisely the
$\Grs$ Hopf algebra. The relation of $\Grs$ with various known
$q$-deformed groups can be exhibited as
\[
\begin{array}{ccccc}
 & & GL_{Q}(3) & &\\
& & \Big\downarrow\vcenter{%
   \rlap{$\mathcal{Q}$}} & & \\
 & & IGL_{Q}(2) & & \\ 
 & & \Big\downarrow\vcenter{%
   \rlap{$\mathcal{Q}$}} & & \\
GL(2)\ot GL(1) &  
\stackrel{\mathcal{L}}{\longleftarrow} &
 \Grs &
\stackrel{\mathcal{F}}{\longrightarrow} &
\GLpqtwo\\
& &\Big\downarrow\vcenter{%
   \rlap{$\mathcal{S}$}} & & \\
 & & \GLqtwo & & 
\end{array}
\] 
where $\mathcal{Q}$, $\mathcal{F}$, $\mathcal{S}$ and $\mathcal{L}$ denote 
the Quotient, Hopf algebra homomorphism, (Hopf)Subalgebra and (classical)
Limit respectively. $GL_{Q}(3)$ denotes the multiparameter $q$-
deformed $GL(3)$ and $IGL_{Q}(2)$ is the inhomogenous multiparameter $q$-
deformation of $GL(2)$. Motivated by the rich structure of $\Grs$, this
quantum group has recently been studied by the authors in detail [2]. As
an intial step in the further understanding of $\Grs$, the authors have
derived explicitly the dual algebra and showed that it is isomorphic to
the single parameter deformation of $gl(2) \oplus gl(1)$, with the second
parameter appearing in the costructure. In [2], the authors have also
constructed a differential calculus on $\Grs$, which inturn provides a
realisation of the calculus on $\GLpqtwo$.
\par

\section{Two parameter $h$- deformations}

Jordanian deformations (also known as $h$-deformations) of Lie groups and
Lie algebras have attracted  a lot of attention in recent years. A
peculiar feature of this deformation is that the corresponding $R$-matrix
is triangular i.e. $R_{12}R_{21}=1$. These deformations are called
`Jordanian' due to the Jordan normal form of the $R$-matrix. It was shown
in [3] that upto isomorphism, $\GLqtwo$ and $\GLhtwo$ are the only
possible distinct deformations (with central determinant) of the group
$GL(2)$. In [4], an interesting observation was made that the $h$-
deformations could be obtained by a singular limit of a similarity
transformation from the $q$- deformations, and this was generalised to
multiparameter deformations as well as to higher dimensions i.e space of
$n \times n$ quantum matrices [5]. For the purpose of current
investigation, the authors have applied the contraction procedure to
$\Grs$ to obtain a new Jordanian quantum group $\Gmk$ [6]. It turns out
that this new structure is also related to other known Jordanian quantum
groups.\\
The $\Gmk$ quantum group can be defined as a triangular Hopf algebra
generated by the $T$- matrix 
$\left( \begin{smallmatrix}a&b&0\\c&d&0\\0&0&f\end{smallmatrix}
\right)$.
The set of commutation relations consisting of elements $a$,$b$,$c$,$d$
form a subalgebra that coincides exactly with the single parameter
Jordanian $\GLhtwo$ for $m=h$. This is exactly analogous to the
$q$-deformed case where the first four elements of $\Grs$ form the
$\GLqtwo$ Hopf subalgebra. Again, the remaining fifth element $f$
generates the $GL(1)$ group, as it did in the $q$-deformed case, and the
second parameter appears only through the cross commutation relations
between $GL_{m}(2)$ and $GL(1)$ elements. Therefore, $\Gmk$ can also be
considered as a two parameter Jordanian deformation of classical
$GL(2)\ot GL(1)$ group. Furthermore, $\Gmk$ also provides a realisation of
the two parameter Jordanian $\GLhh'two$. Besides, it may be
interpreted as a quotient of the multiparameter Jordanian deformation of
$GL(3)$, denoted $GL_{J}(3)$ as well as that of inhomogenous $IGL(2)$,
denoted $IGL_{J}(2)$ quantum groups. This can be represented as follows
\[
\begin{array}{ccccc}
 & & GL_{J}(3) & &\\
& & \Big\downarrow\vcenter{%
   \rlap{$\mathcal{Q}$}} & & \\
 & & IGL_{J}(2) & & \\ 
 & & \Big\downarrow\vcenter{%
   \rlap{$\mathcal{Q}$}} & & \\
GL(2)\ot GL(1) &  
\stackrel{\mathcal{L}}{\longleftarrow} &
 \Gmk &
\stackrel{\mathcal{F}}{\longrightarrow} &
\GLhh'two\\
& &\Big\downarrow\vcenter{%
   \rlap{$\mathcal{S}$}} & & \\
 & & \GLhtwo & & 
\end{array}
\] 
where the maps $\mathcal{Q}$, $\mathcal{F}$, $\mathcal{S}$ and
$\mathcal{L}$ are as before.

\section{Coloured Extensions}

The standard quantum group relations can be extended by parametrising the
corresponding generators using some continuous `colour' variables and
redefining the associated algebra and coalgebra in a way that all Hopf
algebraic properties remain preserved [1,7,8]. For the case of a
single parameter quantum deformation of $\GLtwo$ (with deformation
parameter $r$), its `coloured' version [1] is given by the $R$-matrix,
denoted $R_{r}^{\lambda,\mu}$ which satisfies
\[
R_{12}^{\lambda,\mu}R_{13}^{\lambda,\nu}R_{23}^{\mu,\nu} =
R_{23}^{\mu,\nu}R_{13}^{\lambda,\nu}R_{12}^{\lambda,\mu}
\]
the so-called `Coloured' Quantum Yang Baxter Equation (CQYBE). It should
be stressed at this point that the coloured $R$- matrix provides a
nonadditive-type solution $R^{\lambda,\mu} \neq R(\lambda - \mu)$ of the
Yang-Baxter equation, which is in general
multicomponent and the parameters $\lambda$, $\mu$, $\nu$ are considered
as `colour' parameters. Such solutions were first discovered in the
study of integrable models [9]. This gives rise to the coloured $RTT$
relations 
\[
R_{r}^{\lambda,\mu}T_{1\lambda}T_{2\mu}=T_{2\mu}T_{1\lambda}
R_{r}^{\lambda,\mu}
\]
(where $T_{1\lambda}=T_{\lambda}\dot{\ot} \ident$ and
$T_{2\mu}=\ident\dot{\ot} T_{\mu}$) in which the entries of the $T$
matrices carry colour  dependence. The coproduct and counit for the
coalgebra structure are given by
$\cop (T_{\lambda})=T_{\lambda}\dot{\ot} T_{\lambda}$,
$\cnt (T_{\lambda})=\ident$
and depend only on one colour parameter. By contrast, the algebra 
structure is more complicated with generators of two different 
colours appearing simultaneously in the algebraic relations. The full Hopf
algebraic structure can be constructed and results in a coloured extension
of the quantum group. Since $\lambda$ and $\mu$ are continuous variables,
this implies the coloured quantum group has an infinite number of
generators.\\
The above coloured generalisation of the FRT formalism was given by Kundu
and Basu-Mallick [1,10] and that of the Drinfeld-Jimbo formulation of
quantised universal enveloping algebras has been given by Bonatos, Quesne
{\sl et al} [11]. In the context of knot theory, Ohtsuki [12] introduced
some coloured quasitriangular Hopf algebras, which are characterised by
the existence of a coloured universal $R$- matrix, and he applied his
theory to $U_{q}sl(2)$. Coloured generalisations of quantum groups can
also be understood as an application of the twisting procedure, in a
manner similar to the multiparameter generalisation of quantum groups.
Jordanian deformations also admit coloured extensions  [7]. The associated
$R$-matrix satisfies the CQYBE and is `colour' triangular i.e.
$R_{12}^{\lambda,\mu}=({R_{21}^{\mu,\lambda}})^{-1}$, a coloured extension
of the notion of triangularity.

\subsection*{Coloured Extension of $\Grs$ : $\Grss'$}

The coloured extension of $\Grs$ proposed in [1] has only one deformation
parameter $r$ and two colour paramters $s$ and $s'$. The second
deformation parameter of the uncoloured case now plays the role of a
colour parameter. In such a coloured extension, the first four
generators $a$,$b$,$c$,$d$ are kept independent of the colour parameters
while the fifth generator $f$ is now paramterised by $s$ and $s'$.
The matrices of generators are
\[
T_{s}=\begin{pmatrix}a&b&0\\c&d&0\\0&0&f_{s}\end{pmatrix}\quad , \quad
T_{s'}=\begin{pmatrix}a&b&0\\c&d&0\\0&0&f_{s'}\end{pmatrix}
\]
From the $RTT$ relations, one observes that the commutation relations
between $a$,$b$,$c$,$d$ are as before but $f_{s}$ and $f_{s'}$ now satisfy
two colour copies of the relations satisfied by $f$ of the uncloured
$\Grs$. In addition, the relation $[f_{s},f_{s'}]=0$ holds. The associated
coloured $R$-matrix, denoted $R_{r}^{s,s'}$ satisfies the CQYBE
\[
R_{12}(r;s,s')R_{13}(r;s,s'')R_{23}(r;s',s'') =
R_{23}(r;s',s'')R_{13}(r;s,s'')R_{12}(r;s,s')
\]
and the corresponding coloured quantum group is denoted $\Grss'$.

\subsection*{Coloured Extension of $\Gmk$ : $\Gmkk'$}

Similar to the case of $\Grs$, we have proposed [13] a coloured extension
of the Jordanian quantum group $\Gmk$. The first four generators remain
independent of the colour parameters $k$ and $k'$ whereas the generator
$f$ is parameterised by $k$ and $k'$. Again, the second deformation
parameter $k$ of the uncoloured case now plays the role of a colour
parameter and the $T$-matrices are
\[
T_{k}=\begin{pmatrix}a&b&0\\c&d&0\\0&0&f_{k}\end{pmatrix}\quad , \quad
T_{k'}=\begin{pmatrix}a&b&0\\c&d&0\\0&0&f_{k'}\end{pmatrix}
\]
The commutation relations between $a$,$b$,$c$,$d$ remain unchanged whereas
$f_{k}$ and $f_{k'}$ satisfy two colour copies of the relations satisfied
by $f$ of the uncloured $\Gmk$. In addition, the relation
$[f_{k},f_{k'}]=0$ holds. The associated coloured $R$-matrix,  denoted
$R_{m}^{k,k'}$, is a solution of the CQYBE
\[
R_{12}(m;k,k')R_{13}(m;k,k'')R_{23}(m;k',k'') =
R_{23}(m;k',k'')R_{13}(m;k,k'')R_{12}(m;k,k')
\]
and is colour triangular. The corresponding coloured Jordanian quantum
group is denoted $\Gmkk'$.

\section{Contractions and Homomorphisms}

The $R$-matrix of the Jordanian (or $h$-deformation) can be viewed as a
singular limit of a similarity transformation on the $q$-deformed
$R$-matrix [4]. Let $g(\eta)$ be a matrix dependent on a contraction
parameter $\eta $ which is itself a function of one of the deformation
parameters of the $q$-deformed algebra. This can be used to define a
transformed $q$-deformed $R$-matrix
\[
R_{h} = (g^{-1} \otimes g^{-1})R_{q}(g \otimes g)
\]
The $R$-matrix of the Jordanian deformation is then obtained by taking a
limiting value of the parameter $\eta$. Even though the contraction
parameter $\eta$ is undefined in this limit, the new $R$-matrix is
finite and gives rise to a new quantum group structure through the
$RTT$-relations. For example, in the contraction process which takes
$GL_q(2)$ to $GL_h(2)$, the contraction matrix is 
\[
g(\eta) = \left( \begin{array}{cc} 1 & \eta \\ 0 & 1 \end{array} \right)
\]
where  $\eta = \frac{h}{1-q}$ with $h$ a new free parameter. Such
transformations have proved to be powerful tools in establishing various
connections between the $q$- and the $h$- deformed quantum groups, which
were previously obscure. In the context of the quantum groups under
consideration in the present paper, the contraction procedure was
successfully applied [6] to the $\Grs$ quantum group of Section II to
obtain the Jordanian $\Gmk$ given in Section III. Furthermore, the
multiparameter Jordanian $GL_{J}(3)$ and hence the multiparamter
Inhomogeneous $IGL_{J}(2)$ were also obtained by contracting their
respective $q$- deformed counterparts [14].
\par
The Hopf algebra homomorphism $\mathcal{F}$ from $\Grs$ to $\GLpqtwo$,
which provides a realisation of the latter, is given by
\[
{\mathcal{F}}: \Grs\lmt \GLpqtwo
\]
\[
{\mathcal{F}}\left(\begin{array}{cc}a&b\\c&d\end{array}\right)\lmt
\left(\begin{array}{cc}a'&b'\\c'&d'\end{array}\right)=
f^{N}\left(\begin{array}{cc}a&b\\c&d\end{array}\right)
\]
The elements $a'$,$b'$,$c'$ and $d'$ are the generators of $\GLpqtwo$ and
$N$ is a fixed non-zero integer. The relation between the deformation
parameters $(p,q)$ and $(r,s)$ is given by
\[p = \rinv s^{N} \quad , \quad \quad q = \rinv s^{-N}\]

A Hopf algebra homomorphism 
\[
{\mathcal{F}}: \Gmk\lmt \GLhh'two
\]
of exactly the same form as in the $q$-deformed case, exists between the
generators of $\Gmk$ and $\GLhh'two$ provided that the two sets of
deformation parameters $(h,h')$ and $(m,k)$ are related via the equation
\[
h = -m + Nk \quad , \quad \quad h' = -m - Nk
\]
Note that for vanishing $k$, one gets the one parameter case. In addition,
using the above realisation together with the coproduct, counit and
antipode axioms for the $\Gmk$ algebra and the respective homeomorphism
properties, one can easily recover the standard coproduct, counit and
antipode for $\GLhh'two$. Thus, the Jordanian $\GLhh'two$ group can in
fact be reproduced from the newly defined Jordanian $\Gmk$. It is
curious to note that if we write $p=e^{h}$, $q=e^{h'}$, $r=e^{m}$ and
$s=e^{k}$, then the relations between the parameters in the $q$-deformed
case and the $h$-deformed case are identical. The systematics of the
uncloured quantum groups discussed here can be summarised in the following
commutative diagram
\[
\begin{CD}
GL_{Q}(3) @>\mathcal{Q}>> IGL_{Q}(2) @>\mathcal{Q}>> \Grs
@>\mathcal{F}>> \GLpqtwo\\
@V{\mathcal{C}}VV    @V{\mathcal{C}}VV @VV{\mathcal{C}}V @VV{\mathcal{C}}V
\\
GL_{J}(3) @>>\mathcal{Q}> IGL_{J}(2) @>>\mathcal{Q}> \Gmk @>>\mathcal{F}>
\GLhh'two
\end{CD}
\]
where $\mathcal{Q}$, $\mathcal{C}$ and $\mathcal{F}$ denote the
quotient, contraction and the Hopf algebra homomorphism. The contraction
procedure discussed above has been successfully applied [13] to the case
of coloured quantum groups yielding new coloured Jordanian deformations.
We apply to $R_{r}^{\lambda,\mu}$, the coloured $R$-matrix for
$q$-deformed $\GLtwo$, the transformation
\[
(g\ot g)^{-1}R_{r}^{\lambda,\mu}(g\ot g)
\]
where $g$ is the two dimensional transformation matrix
$\left( \begin{smallmatrix}1&\eta\\0&1\end{smallmatrix} \right)$ and
$\eta$ is chosen to be $\eta=\frac{m}{1-r}$. In the limit $r\rightarrow
1$, we obtain a new $R$-matrix, $R_{m}^{\lambda,\mu}$ which is a coloured
$R$-matrix for a Jordanian deformation of $\GLtwo$. The contraction is
then also used to obtain the coloured extension $\Gmkk'$ of $\Gmk$, from
the coloured extension $\Grss'$ of $\Grs$.  The $R$-matrix $R_{m}^{k,k'}$
is obtained as the contraction limit of the $R$-matrix for the coloured
extension of $\Grs$ via the transformation
\[
R_{m}^{k,k'} = \lim_{r\rightarrow 1}(G\ot G)^{-1} R_{r}^{s,s'}(G\ot G)
\]
where
\[
G=\begin{pmatrix}g&0\\0&1\end{pmatrix};\quad
g=\begin{pmatrix}1&\eta\\0&1\end{pmatrix},\quad \eta=\frac{m}{r-1}
\]
The Hopf algebra homomorphism from $\Grss'$ to $GL_{r}^{\lambda,\mu}(2)$
\[
\mathcal{F}_{N}: \Grss'\lmt GL_{r}^{\lambda,\mu}(2)
\]
is given by
\[
\mathcal{F}_{N}:\begin{pmatrix}a&b\\c&d\end{pmatrix}\lmt
\begin{pmatrix}a'_{\lambda}&b'_{\lambda}\\c'_{\lambda}&d'_{\lambda}
\end{pmatrix}=f_{s}^{N}\begin{pmatrix}a&b\\c&d\end{pmatrix}
\]
\[
\mathcal{F}_{N}:\begin{pmatrix}a&b\\c&d\end{pmatrix}\lmt
\begin{pmatrix}a'_{\mu}&b'_{\mu}\\c'_{\mu}&d'_{\mu}
\end{pmatrix}=f_{s'}^{N}\begin{pmatrix}a&b\\c&d\end{pmatrix}
\]
where N is a fixed non-zero integer and the sets of colour parameters
$(s,s')$ and $(\lambda,\mu)$ are related through quantum deformation
parameter $r$ by
\[
s=r^{2N\lambda} \quad , \quad s'=r^{2N\mu}
\]
The primed generators $a'_{\lambda}$,
$b'_{\lambda}$,$c'_{\lambda}$,$d'_{\lambda}$ and
$a'_{\mu}$,$b'_{\mu}$,$c'_{\mu}$,$d'_{\mu}$ belong to
$GL_{r}^{\lambda,\mu}(2)$ whereas the unprimed ones
$a$,$b$,$c$,$d$,$f_{s}$ and  $f_{s'}$ are generators of $\Grss'$.
If we now denote the generators of   
$GL_{m}^{\lambda,\mu}(2)$ by
$a'_{\lambda}$,$b'_{\lambda}$,$c'_{\lambda}$,$d'_{\lambda}$ and 
$a'_{\mu}$,$b'_{\mu}$,$c'_{\mu}$,$d'_{\mu}$ and the generators of $\Gmkk'$
by $a$,$b$,$c$,$d$,$f_{k}$ and  $f_{k'}$ then a Hopf algebra homomorphism
from $\Gmkk'$ to
$GL_{m}^{\lambda,\mu}(2)$
\[
\mathcal{F}_{N}: \Gmkk'\lmt GL_{m}^{\lambda,\mu}(2)
\]
is of exactly the same form
\[
\mathcal{F}_{N}:\begin{pmatrix}a&b\\c&d\end{pmatrix}\lmt
\begin{pmatrix}a'_{\lambda}&b'_{\lambda}\\c'_{\lambda}&d'_{\lambda}
\end{pmatrix}=f_{k}^{N}\begin{pmatrix}a&b\\c&d\end{pmatrix}
\]   
\[
\mathcal{F}_{N}:\begin{pmatrix}a&b\\c&d\end{pmatrix}\lmt
\begin{pmatrix}a'_{\mu}&b'_{\mu}\\c'_{\mu}&d'_{\mu}
\end{pmatrix}=f_{k'}^{N}\begin{pmatrix}a&b\\c&d\end{pmatrix}
\]
The sets of colour parameters $(k,k')$ and $(\lambda,\mu)$ are related
to the Jordanian deformation parameter $m$ by
\[
Nk=-2m\lambda \quad , \quad Nk'=-2m\mu
\]   
and N, again, is a fixed non-zero integer. The schematics of our analysis
for the coloured quantum groups is represented in the diagram
\[
\begin{CD}
\Grs @>\mathcal{E}>> \Grss' @>\mathcal{F}>> GL_{r}^{\lambda,\mu}(2)\\
@V{\mathcal{C}}VV    @VV{\mathcal{C}}V @VV{\mathcal{C}}V\\
\Gmk @>>\mathcal{E}> \Gmkk' @>>\mathcal{F}> GL_{m}^{\lambda,\mu}(2)
\end{CD}
\]
where $\mathcal{C}$, $\mathcal{F}$ and $\mathcal{E}$ denote the
contraction, Hopf algebra homomorphism and coloured extension 
respectively. In both of the commutative diagrams above, the objects at
the top level are the $q$ deformed ones and the corresponding Jordanian
counterparts are shown at the bottom level.

\section{Conclusions}

In the present work, we have obtained a new Jordanian quantum group
$G_{m,k}$ by contraction of the $q$- deformed quantum group $G_{r,s}$. We
then used this new structure to establish quantum group homomorphisms with
other known two parameter quantum groups at the Jordanian level. At the
same time we also showed that such homomorphisms commute with the
contraction procedure. Our analysis is then set in the wider context
of coloured quantum groups. We give a coloured generalisation of the
contraction procedure and obtain new coloured Jordanian quantum groups.
A careful study of the properties of both $G_{r,s}$ and $G_{m,k}$ lead to
their respective coloured extensions. Furthermore, we show that the
homomorphisms of the uncoloured case naturally extend to the coloured
case.
\par
The physical interest in studying $\Grs$ lies in the observation that when
endowed with a $\ast$- structure, this quantum group specialises to a
two parameter quantum deformation of $SU(2) \ot U(1)$ which is precisely
the gauge group for the theory of electroweak interactions. Since gauge
theories have an obvious differential geometric description, the study of
differential calculus [2] provides insights in constructing a $q$-gauge
theory based on $\Grs$. It would also be of significance to generalise the
formalism of differential calculus to the case of coloured quantum
groups and explore possible physical applications.
\par

\section*{Acknowledgments}

D.P. is grateful to the organisers of the Symposium and would like to
thank Prof. David Radford and Dr. Gustav Delius for useful comments. The
authors have also benefited by discussions with Prof. Vlado Dobrev and
Dr. Preeti Parashar.\par

\section*{References}

[1] B. Basu-Mallick, hep-th/9402142; {\sl Intl. J. Mod. Phys.} {\bf A10},
2851 (1995).\par
[2] D. Parashar and R. J. McDermott, Kyoto University preprint RIMS - 1260
(1999), math.QA/9901132.\par
[3] B. A. Kupershmidt, {\sl J. Phys.} {\bf A25}, L1239 (1992).\par
[4] A. Aghamohammadi, M. Khorrami and A. Shariati, {\sl J. Phys.} {\bf
A28}, L225 (1995).\par
[5] M. Alishahiha, {\sl J. Phys.} {\bf A28}, 6187 (1995).\par
[6] D. Parashar and R. J. McDermott, math.QA/9909001, {\sl Czech. J.
Phys.}, in press.\par
[7] P. Parashar, {\sl Lett. Math. Phys.} {\bf 45}, 105 (1998).\par
[8] C. Quesne, {\sl J. Math. Phys.} {\bf 38}, 6018 (1997); {\em ibid} {\bf 
39}, 1199 (1998).\par
[9] V. V. Bazhanov and Yu. G. Stroganov, {\sl Theor. Math. Phys.} {\bf
62}, 253 (1985).\par
[10] A. Kundu and B. Basu-Mallick, {\sl J. Phys.} {\bf A27}, 3091 (1994);
B. Basu-Mallick, {\sl Mod. Phys. Lett.} {\bf A9}, 2733 (1994).\par
[11] D. Bonatos {\sl et. al.}, {\sl J. Math. Phys.} {\bf 38}, 369 (1997);
C. Quesne, q-alg/9705022.\par
[12] T. Ohtsuki, {\sl J. Knot Theor. Its Rami.} {\bf 2}, 211 (1993).\par
[13] D. Parashar and R. J. McDermott, math.QA/9911194, {\sl J. Math.
Phys.}, in press.\par
[14] R. J. McDermott and D. Parashar, math.QA/9909045, {\sl Czech. J.
Phys.}, in press.\par

\end{document}